\documentclass[11pt]{article}
\usepackage{fullpage}
\usepackage{amsmath}
\usepackage{amsfonts}
\usepackage{amsthm}
\usepackage{graphicx}
\usepackage{caption}
\usepackage[table,xcdraw]{xcolor}
\usepackage{booktabs}
\usepackage{mdframed}
\usepackage{subcaption}
\usepackage{enumerate}
\usepackage{hyperref}
\usepackage{mathtools}
\usepackage{comment}
\usepackage{tikz}
\usepackage{float}
\usepackage{lastpage}
\usepackage{enumerate}
\usepackage{fancyhdr}
\usepackage{mathrsfs}
\usepackage{xcolor}
\usepackage{listings}
\usepackage{verbatim}
\usepackage{bm}
\usepackage{extarrows}  % \xlongequal[\quad\quad]{u=r^2}
\usetikzlibrary{positioning}
\usetikzlibrary{arrows}

\usepackage[parfill]{parskip} %so paragraphs don't indent

\usepackage{fancyhdr}
\numberwithin{equation}{section}

\begin{document}

\title{A Mathematical Analysis of Mathematical Faculty}
\author{Victoria Chayes \and Dodam Ih \and Yukun Yao \and Doron Zeilberger \and Tianhao Zhang}

\maketitle

\leavevmode
\\

\begin{abstract}
 We use the data of tenured and tenure-track faculty at ten public and private math departments of various tiered rankings in the United States, as a case study to demonstrate the statistical and mathematical relationships among several variables, e.g., the number of publications and citations, the rank of professorship and AMS fellow status. At first we do an exploratory data analysis of the math departments. Then various statistical tools, including regression, artificial neural network, and unsupervised learning, are applied and the results obtained from different methods are compared. We conclude that with more advanced models, it may be possible to design an automatic promotion algorithm that has the potential to be fairer, more efficient and more consistent than human approach. 

\end{abstract}
\leavevmode
\\
\\

\section{Introduction}
Modern research universities and colleges around the globe employ tenure track and tenured professors in STEM fields in a large part for the quality and impact of the research they produce. However, the process of promotion on this tenure track is inefficient, time-consuming, and varies from institution to institution. It is largely based on subjective evaluations by ``experts", and endless committee meetings, with many steps along the way where personal bias may undermine a promising candidate. As such, an algorithmic approach to promotion may be a vast improvement to the current system. The purpose of this paper is to determine whether or not this is feasible, by testing if automated algorithms can predict the rank of tenure track professors in public and private universities in the United States. As the authors are from a department of mathematics and are most familiar with mathematical academia, we concentrate on math departments.

The variables considered in the development of a promotion algorithm for mathematics are: years from PhD, number of papers, number of citations, h-index, and AMS Fellowship status. We attempt to predict if a candidate is an associate professor, assistant professor, professor, or distinguished professor from this, using tools of regression, artificial neural network and unsupervised learning to the data set we collected for predictive exploration. We discover (perhaps unsurprisingly) that there are very strong statistical properties shared within certain groups of ranks. We conclude this means that an unbiased algorithm could be built by studying data and finding patterns to avoid subjective opinions and maintain consistent standards. It is important to note that while race and gender of the candidates was beyond the scope of this study, instituting an algorithmic approach to promotion may also help mitigate negative bias with regards to race or gender in academic promotions.

For this paper, we collect public data online of the tenured and tenure-track faculty members at math departments of 10 universities: UC Berkeley, Dartmouth College, University of Florida, Harvard University, University of Michigan, Massachusetts Institute of Technology, University of Pennsylvania, Princeton University, Rutgers University-New Brunswick and UCLA.  In our data set, totally there are 444 professors. The variables are Lastname, Firstname, Rank (Assistant Professor=1, Associate Professor=2, Full Professor=3, Distinguished Professor=4), Number of Publication, Number of Citations, $h$-Index, AMS Fellow (Yes=1, No=0), the Year of Ph.D. Awarded. 

The data of names and ranks are from the website of each math department. Note that some departments, e.g., Princeton and Harvard, do not have the rank of Distinguished Professor. The numbers of publications, citations and $h$-index are from MathSciNet. It is worth emphasizing that MathSciNet has much more strict standards for recording publications, citations and $h$-index so that the MathSciNet $h$-index is much lower (roughly a half) from that which can be found using Google Scholar. AMS membership was taken from the AMS website. Year of attaining the rank of PhD was taken from Mathematics Genealogy. All data was collected in or around November 2019.

\section{Exploratory Data Analysis}

The following three charts show the means and the standard deviations for each of the fields considered across all universities.

\begin{table}[h!]
	\centering
	\begin{tabular}{|r|c | c|}
		
		\hline	
		Field & Mean & Standard Deviation \\
		\hline
		
		Rank & 2.757 & 0.763\\
		Number of Publications & 62.459 & 60.063\\
		Number of Citations & 1250.153 & 2012.400\\
		\textit{h}-index & 14.777 & 9.866\\
		AMS Fellowship & 0.358 & 0.480\\
		Year of PhD & 1992.304 & 14.538\\
		\hline	
	\end{tabular}
	\caption{Means and standard deviations across all universities (n=444)}
\end{table}

\begin{table}[h!]
	\small
	\centering
	\begin{tabular}{|r |c |c c c c c c|}
		
		\hline	
University & n & Rank & Publications & Citations & \textit{h}-index & AMS Fellowship & Year of PhD \\

		\hline
		
Berkeley & 58 & 2.741 & 64.914 & 1579.017 & 17.207 & 0.362 & 1992.776 \\ \hline
Dartmouth & 23 & 2.478 & 36.783 & 360.435 & 8.652 & 0.043 & 1993.652 \\ \hline
Florida & 44 & 2.500 & 50.568 & 416.477 & 9.136 & 0.045 & 1992.091 \\ \hline
Harvard & 20 & 3.000 & 100.400 & 2810.800 & 24.500 & 0.400 & 1984.000 \\ \hline
MIT & 53 & 2.642 & 63.491 & 1460.094 & 16.377 & 0.415 & 1995.717 \\ \hline
Michigan & 62 & 2.871 & 54.258 & 936.742 & 12.887 & 0.339 & 1991.694 \\ \hline
Penn & 25 & 2.800 & 53.960 & 633.200 & 12.320 & 0.400 & 1989.440 \\ \hline
Princeton & 42 & 2.452 & 73.524 & 2123.738 & 19.357 & 0.452 & 1995.071 \\ \hline
Rutgers & 59 & 3.153 & 71.661 & 1027.525 & 14.271 & 0.559 & 1989.051 \\ \hline
UCLA & 58 & 2.776 & 60.241 & 1371.379 & 14.517 & 0.379 & 1994.397 \\ \hline

	\end{tabular}
	\caption{Means across all universities, by university (n=444)}
\end{table}

\begin{table}[h!]
	\small
	\centering
	\begin{tabular}{|r |c |c c c c c c|}
		
		\hline	
	University & n & Rank & Publications & Citations & \textit{h}-index & AMS Fellowship & Year of PhD\\
		
		\hline
	Berkeley & 58 & 0.609 & 48.665 & 2174.119 & 9.472 & 0.485 & 12.445 \\ \hline
	Dartmouth & 23 & 0.790 & 30.705 & 399.379 & 4.914 & 0.209 & 12.463 \\ \hline
	Florida & 44 & 0.876 & 37.279 & 507.301 & 5.129 & 0.211 & 15.397 \\ \hline
	Harvard & 20 & 0.000 & 106.460 & 3291.795 & 11.390 & 0.503 & 12.645 \\ \hline
	MIT & 53 & 0.787 & 59.765 & 2083.936 & 10.895 & 0.497 & 15.468 \\ \hline
	Michigan & 62 & 0.614 & 45.168 & 1364.324 & 7.378 & 0.477 & 13.745 \\ \hline
	Penn & 25 & 0.645 & 31.798 & 465.785 & 5.429 & 0.500 & 14.509 \\ \hline
	Princeton & 42 & 0.889 & 90.541 & 2465.433 & 13.483 & 0.504 & 17.374 \\ \hline
	Rutgers & 59 & 0.979 & 62.757 & 1128.886 & 7.850 & 0.501 & 15.234 \\ \hline
	UCLA & 58 & 0.531 & 59.434 & 2898.606 & 10.881 & 0.489 & 13.124 \\ \hline

	\end{tabular}
	\caption{Standard deviations across all universities, by university (n=444)}
\end{table}

Out of the mathematics departments in the study, Harvard led in the averages for number of publications, number of citations, the \textit{h}-index, and number of years since Ph.D., followed in each of these metrics except for the last by Princeton, whose faculty on average completed their doctorates a full eleven years after their Harvard colleagues. Princeton also had the greatest variation in the $h$-index and the academic age of its faculty; however, this may result from a different classification system used by Princeton that does not award full tenure. Rutgers led in both the average rank of the faculty titles and the proportion of AMS Fellows, despite being slightly below average in both the number of citations and the \textit{h}-index.

Dartmouth, Florida, and Penn had the lowest variability among its faculty in the number of publications, number of citations, and the \textit{h}-index.

The following charts show the extrapolated percentiles for each field:
\begin{table}[h!]
	
	\centering
	\begin{tabular}{|r|c c c c c c c|}
		
		\cline{2-8}
		\multicolumn{1}{c|}{} & \multicolumn{7}{c|}{Percentiles} \\
		\hline
		Field & \multicolumn{1}{c}{5} &  \multicolumn{1}{|c}{10} & \multicolumn{1}{|c}{25} & \multicolumn{1}{|c}{50} & \multicolumn{1}{|c}{75} & \multicolumn{1}{|c}{90} & \multicolumn{1}{|c|}{95} \\
		\hline
		
		Rank & 2 & 2 & 3 & 3 & 3 & 3 & 3 \\ \hline
		Publications &  11 & 14 & 25 & 44 & 71 & 111 & 164 \\ \hline
		Citations & 30 & 68 & 208 & 626 & 1392 & 2260 & 3750 \\ \hline
		\textit{h}-index  & 4 & 4 & 8 & 13 & 18 & 22 & 31 \\ \hline
		AMS Fellowship  & 0 & 0 & 0 & 0 & 1 & 1 & 1 \\ \hline
		Year of PhD &  1970 & 1977 & 1986 & 1996 & 2006 & 2009 & 2012\\
		\hline	
	\end{tabular}
	\caption{Percentiles for all fields across all universities (n=444)}
\end{table}

\begin{table}[h!]
	
	\centering
	\begin{tabular}{|r|c c c c c c c|}
		
		\cline{2-8}
		\multicolumn{1}{c|}{} & \multicolumn{7}{c|}{Percentiles} \\
		\hline
		University & \multicolumn{1}{c}{5} &  \multicolumn{1}{|c}{10} & \multicolumn{1}{|c}{25} & \multicolumn{1}{|c}{50} & \multicolumn{1}{|c}{75} & \multicolumn{1}{|c}{90} & \multicolumn{1}{|c|}{95} \\
		\hline
Berkeley & 1 & 2 & 3 & 3 & 3 & 3 & 3\\ \hline
Dartmouth & 1 & 1 & 2 & 3 & 3 & 3 & 3\\ \hline
Florida & 1 & 1 & 2 & 3 & 3 & 3 & 4\\ \hline
Harvard & 3 & 3 & 3 & 3 & 3 & 3 & 3\\ \hline
MIT & 1 & 1 & 3 & 3 & 3 & 3 & 3\\ \hline
Michigan & 2 & 2 & 3 & 3 & 3 & 3 & 4\\ \hline
Penn & 1 & 2 & 3 & 3 & 3 & 3 & 3\\ \hline
Princeton & 1 & 1 & 1 & 3 & 3 & 3 & 3\\ \hline
Rutgers & 1 & 2 & 3 & 3 & 4 & 4 & 4\\ \hline
UCLA & 2 & 2 & 3 & 3 & 3 & 3 & 3\\ \hline
	\end{tabular}

	\caption{Percentiles for professor rank across all universities, by university (n=444)}
\end{table}

\begin{table}[h!]
	
	\centering
	\begin{tabular}{|r|c c c c c c c|}
		
		\cline{2-8}
		\multicolumn{1}{c|}{} & \multicolumn{7}{c|}{Percentiles} \\
		\hline
		University & \multicolumn{1}{c}{5} &  \multicolumn{1}{|c}{10} & \multicolumn{1}{|c}{25} & \multicolumn{1}{|c}{50} & \multicolumn{1}{|c}{75} & \multicolumn{1}{|c}{90} & \multicolumn{1}{|c|}{95} \\
		\hline
		Berkeley & 18 & 23 & 30 & 43 & 88 & 131 & 139\\ \hline
		Dartmouth & 9 & 9 & 14 & 26 & 58 & 65 & 75\\ \hline
		Florida & 8 & 11 & 19 & 42 & 76 & 98 & 106\\ \hline
		Harvard & 38 & 39 & 48 & 74 & 112 & 137 & 158\\ \hline
		MIT & 13 & 14 & 26 & 40 & 81 & 133 & 215\\ \hline
		Michigan & 10 & 13 & 21 & 45 & 70 & 100 & 122\\ \hline
		Penn & 11 & 18 & 28 & 54 & 75 & 103 & 110\\ \hline
		Princeton & 3 & 6 & 10 & 53 & 94 & 160 & 207\\ \hline
		Rutgers & 10 & 15 & 26 & 58 & 104 & 138 & 155\\ \hline
		UCLA & 11 & 14 & 25 & 44 & 71 & 111 & 164\\ \hline
	\end{tabular}
	
	\caption{Percentiles for number of publications across all universities, by university (n=444)}
\end{table}

\begin{table}[h!]
	
	\centering
	\begin{tabular}{|r|c c c c c c c|}
		
		\cline{2-8}
		\multicolumn{1}{c|}{} & \multicolumn{7}{c|}{Percentiles} \\
		\hline
		University & \multicolumn{1}{c}{5} &  \multicolumn{1}{|c}{10} & \multicolumn{1}{|c}{25} & \multicolumn{1}{|c}{50} & \multicolumn{1}{|c}{75} & \multicolumn{1}{|c}{90} & \multicolumn{1}{|c|}{95} \\
		\hline
		
		Berkeley & 188 & 210 & 366 & 718 & 1786 & 3800 & 5320\\ \hline
		Dartmouth & 35 & 38 & 72 & 167 & 510 & 918 & 1086\\ \hline
		Florida & 30 & 33 & 70 & 261 & 488 & 965 & 1361\\ \hline
		Harvard & 743 & 754 & 990 & 1440 & 3000 & 4862 & 9490\\ \hline
		MIT & 28 & 52 & 165 & 622 & 1683 & 3194 & 4968\\ \hline
		Michigan & 31 & 76 & 184 & 602 & 1087 & 1759 & 3501\\ \hline
		Penn & 56 & 99 & 241 & 586 & 927 & 1285 & 1379\\ \hline
		Princeton & 9 & 37 & 142 & 1158 & 3664 & 5732 & 6260\\ \hline
		Rutgers & 38 & 75 & 209 & 632 & 1446 & 2536 & 3500\\ \hline
		UCLA & 30 & 68 & 208 & 626 & 1392 & 2260 & 3750\\ \hline
	\end{tabular}
	
	\caption{Percentiles for number of citations across all universities, by university (n=444)}
\end{table}

\begin{table}[h!]
	
	\centering
	\begin{tabular}{|r|c c c c c c c|}
		
		\cline{2-8}
		\multicolumn{1}{c|}{} & \multicolumn{7}{c|}{Percentiles} \\
		\hline
		University & \multicolumn{1}{c}{5} &  \multicolumn{1}{|c}{10} & \multicolumn{1}{|c}{25} & \multicolumn{1}{|c}{50} & \multicolumn{1}{|c}{75} & \multicolumn{1}{|c}{90} & \multicolumn{1}{|c|}{95} \\
		\hline
Berkeley & 7 & 8 & 9 & 15 & 22 & 33 & 36\\ \hline
Dartmouth & 3 & 3 & 5 & 8 & 11 & 16 & 18\\ \hline
Florida & 3 & 3 & 4 & 9 & 12 & 16 & 17\\ \hline
Harvard & 14 & 14 & 17 & 22 & 32 & 35 & 39\\ \hline
MIT & 3 & 6 & 9 & 13 & 23 & 29 & 34\\ \hline
Michigan & 3 & 4 & 8 & 12 & 17 & 23 & 26\\ \hline
Penn & 5 & 5 & 8 & 12 & 16 & 20 & 20\\ \hline
Princeton & 2 & 4 & 6 & 16 & 31 & 39 & 41\\ \hline
Rutgers & 4 & 5 & 8 & 13 & 18 & 26 & 27\\ \hline
UCLA & 4 & 4 & 8 & 13 & 18 & 22 & 31\\ \hline
	\end{tabular}
	
	\caption{Percentiles for the \textit{h}-index across all universities, by university (n=444)}
\end{table}

\begin{table}[h!]
	
	\centering
	\begin{tabular}{|r|c c c c c c c|}
		
		\cline{2-8}
		\multicolumn{1}{c|}{} & \multicolumn{7}{c|}{Percentiles} \\
		\hline
		University & \multicolumn{1}{c}{5} &  \multicolumn{1}{|c}{10} & \multicolumn{1}{|c}{25} & \multicolumn{1}{|c}{50} & \multicolumn{1}{|c}{75} & \multicolumn{1}{|c}{90} & \multicolumn{1}{|c|}{95} \\
		\hline
Berkeley & 1973 & 1975 & 1984 & 1992 & 2002 & 2010 & 2011\\ \hline
Dartmouth & 1979 & 1979 & 1982 & 1996 & 2004 & 2010 & 2012\\ \hline
Florida & 1969 & 1973 & 1982 & 1988 & 2006 & 2013 & 2014\\ \hline
Harvard & 1966 & 1967 & 1978 & 1984 & 1991 & 2000 & 2004\\ \hline
MIT & 1966 & 1974 & 1988 & 1997 & 2008 & 2013 & 2014\\ \hline
Michigan & 1967 & 1971 & 1983 & 1994 & 2002 & 2009 & 2011\\ \hline
Penn & 1965 & 1976 & 1980 & 1987 & 2002 & 2009 & 2012\\ \hline
Princeton & 1966 & 1973 & 1980 & 2000 & 2011 & 2014 & 2015\\ \hline
Rutgers & 1968 & 1969 & 1977 & 1989 & 2000 & 2011 & 2014\\ \hline
UCLA & 1970 & 1977 & 1986 & 1996 & 2006 & 2009 & 2012\\ \hline
	\end{tabular}
	
	\caption{Percentiles for year of PhD across all universities, by university (n=444)}
\end{table}

The covariance and correlation matrices follow:

\begin{table}[h!]
	\centering
	\begin{tabular}{ |c|c|c|c|c|c|c|}
		\hline	
		& Rank  & Publications & Citations & \textit{h}-index & AMS Fellowship & Year of PhD   \\
		\hline
	Rank & 0.582 & 19.624 & 437.279 & 3.485 & 0.151 & -7.125 \\ \hline
	Publications & 19.624 & 3607.558 & 94206.857 & 473.972 & 10.600 & -461.266 \\ \hline
	Citations & 437.279 & 94206.857 & 4049753.638 & 17365.165 & 352.855 & -12455.609 \\ \hline
	\textit{h}-index & 3.485 & 473.972 & 17365.165 & 97.338 & 2.166 & -73.794 \\ \hline
	AMS Fellowship & 0.151 & 10.600 & 352.855 & 2.166 & 0.230 & -2.617 \\ \hline
	Year of PhD & -7.125 & -461.266 & -12455.609 & -73.794 & -2.617 & 211.359 \\ \hline
	
	\end{tabular}
	\caption{The covariance matrix for all universities (n=444)}
\end{table}

\begin{table}[h!]
	\centering
	\begin{tabular}{ |c|c|c|c|c|c|c|}
		\hline	
		& Rank  & Publications & Citations & \textit{h}-index & AMS Fellowship & Year of PhD   \\
		\hline
	Rank & 1.000 & 0.428 & 0.285 & 0.463 & 0.411 & -0.642 \\ \hline
	Publications & 0.428 & 1.000 & 0.779 & 0.800 & 0.368 & -0.528 \\ \hline
	Citations & 0.285 & 0.779 & 1.000 & 0.875 & 0.365 & -0.426 \\ \hline
	\textit{h}-index & 0.463 & 0.800 & 0.875 & 1.000 & 0.457 & -0.514 \\ \hline
	AMS Fellowship & 0.411 & 0.368 & 0.365 & 0.457 & 1.000 & -0.375 \\ \hline
	Year of PhD & -0.642 & -0.528 & -0.426 & -0.514 & -0.375 & 1.000 \\ \hline
		
	\end{tabular}
	\caption{The correlation matrix for all universities (n=444)}
\end{table}

Dividing the universities surveyed into two groups depending on whether they are private or public, we have the following comparisons:

\begin{table}[h!]
	\centering
	\begin{tabular}{|r|c | c|}
		\hline	
		Field & Mean & Standard Deviation \\
		\hline
		Rank & 2.638 & 0.760\\
		Number of Publications & 65.374 & 71.654\\
		Number of Citations & 1514.834 & 2206.948\\
		\textit{h}-index & 16.429 & 11.333\\
		AMS Fellowship & 0.368 & 0.484\\
		Year of PhD & 1992.859 & 15.484\\
		\hline	
	\end{tabular}
	\caption{Means and standard deviations across private universities (n=163)}
\end{table}

\begin{table}[h!]
	\centering
	\begin{tabular}{|r|c | c|}
		
		\hline	
		Field & Mean & Standard Deviation \\
		\hline
		
		Rank & 2.826 & 0.757\\
		Number of Publications & 60.769 & 52.242\\
		Number of Citations & 1096.619 & 1877.458\\
		\textit{h}-index & 13.819 & 8.785\\
		AMS Fellowship & 0.352 & 0.479\\
		Year of PhD & 1991.982 & 13.979\\
		\hline	
	\end{tabular}
	\caption{Means and standard deviations across public universities (n=281)}
\end{table}

\begin{table}[h!]
	\centering
	\begin{tabular}{ |c|c|c|c|c|c|c|}
		\hline	
		& Rank  & Publications & Citations & \textit{h}-index & AMS Fellowship & Year of PhD   \\
		\hline
		Rank & 0.578 & 19.815 & 535.983 & 4.200 & 0.146 & -7.379 \\ \hline
		Publications & 19.815 & 5134.347 & 139945.667 & 664.980 & 12.923 & -586.052 \\ \hline
		Citations & 535.983 & 139945.667 & 4870620.497 & 22485.337 & 377.932 & -17699.863 \\ \hline
		\textit{h}-index & 4.200 & 664.980 & 22485.337 & 128.432 & 2.390 & -103.130 \\ \hline
		AMS Fellowship & 0.146 & 12.923 & 377.932 & 2.390 & 0.234 & -3.090 \\ \hline
		Year of PhD & -7.379 & -586.052 & -17699.863 & -103.130 & -3.090 & 239.752 \\ \hline
		
	\end{tabular}
	\caption{The covariance matrix for private universities (n=163)}
\end{table}

\begin{table}[h!]
	\centering
	\begin{tabular}{ |c|c|c|c|c|c|c|}
		\hline	
		& Rank  & Publications & Citations & \textit{h}-index & AMS Fellowship & Year of PhD   \\
		\hline
		Rank & 0.573 & 19.902 & 410.637 & 3.265 & 0.155 & -6.942 \\ \hline
		Publications & 19.902 & 2729.271 & 67370.508 & 360.722 & 9.268 & -392.204 \\ \hline
		Citations & 410.637 & 67370.508 & 3524847.394 & 14062.499 & 337.174 & -9601.000 \\ \hline
		\textit{h}-index & 3.265 & 360.722 & 14062.499 & 77.185 & 2.028 & -57.928 \\ \hline
		AMS Fellowship & 0.155 & 9.268 & 337.174 & 2.028 & 0.229 & -2.358 \\ \hline
		Year of PhD & -6.942 & -392.204 & -9601.000 & -57.928 & -2.358 & 195.403 \\ \hline
	\end{tabular}
	\caption{The covariance matrix for public universities (n=281)}
\end{table}

\begin{table}[h!]
	\centering
	\begin{tabular}{ |c|c|c|c|c|c|c|}
		\hline	
		& Rank  & Publications & Citations & \textit{h}-index & AMS Fellowship & Year of PhD   \\
		\hline
Rank & 1.000 & 0.364 & 0.319 & 0.487 & 0.398 & -0.627 \\ \hline
Publications & 0.364 & 1.000 & 0.885 & 0.819 & 0.373 & -0.528 \\ \hline
Citations & 0.319 & 0.885 & 1.000 & 0.899 & 0.354 & -0.518 \\ \hline
\textit{h}-index & 0.487 & 0.819 & 0.899 & 1.000 & 0.436 & -0.588 \\ \hline
AMS Fellowship & 0.398 & 0.373 & 0.354 & 0.436 & 1.000 & -0.412 \\ \hline
Year of PhD & -0.627 & -0.528 & -0.518 & -0.588 & -0.412 & 1.000 \\ \hline
	\end{tabular}
	\caption{The correlation matrix for private universities (n=163)}
\end{table}

\begin{table}[h!]
	\centering
	\begin{tabular}{ |c|c|c|c|c|c|c|}
		\hline	
		& Rank  & Publications & Citations & \textit{h}-index & AMS Fellowship & Year of PhD   \\
		\hline
Rank & 1.000 & 0.503 & 0.289 & 0.491 & 0.427 & -0.656 \\ \hline
Publications & 0.503 & 1.000 & 0.687 & 0.786 & 0.371 & -0.537 \\ \hline
Citations & 0.289 & 0.687 & 1.000 & 0.853 & 0.375 & -0.366 \\ \hline
\textit{h}-index & 0.491 & 0.786 & 0.853 & 1.000 & 0.482 & -0.472 \\ \hline
AMS Fellowship & 0.427 & 0.371 & 0.375 & 0.482 & 1.000 & -0.352 \\ \hline
Year of PhD & -0.656 & -0.537 & -0.366 & -0.472 & -0.352 & 1.000 \\ \hline
	\end{tabular}
	\caption{The correlation matrix for public universities (n=281)}
\end{table}

The correlation matrices were largely similar between the private and public universities studied, with three notable exceptions:
\begin{itemize}
	\item  The correlation between the faculty rank and the number of publications was much stronger for public universities (0.503) than for private universities (0.364).
	\item The correlation between the number of publications and the number of citations was much stronger for private universities (0.885) than for public universities (0.687).
	\item The correlation between the academic age and the number of citations was much stronger for private universities (-0.518) than for public universities (-0.366).
\end{itemize}

\begin{figure}[h!]
	\centering
	\includegraphics[width=5cm]{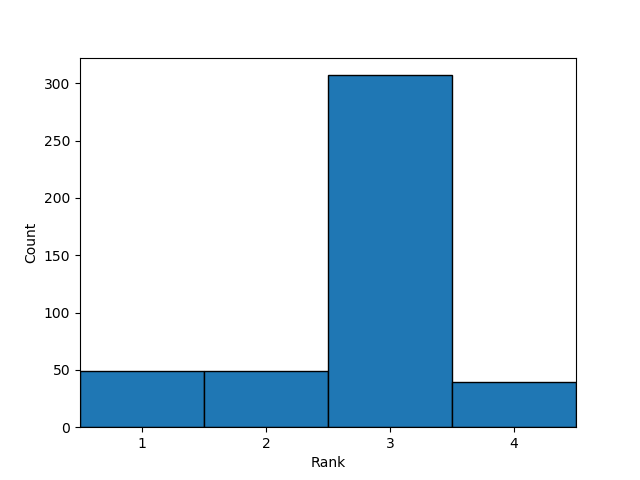}
	\includegraphics[width=5cm]{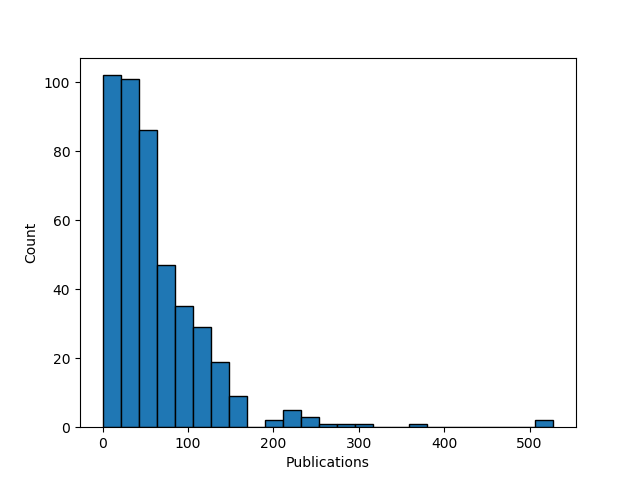}
	\includegraphics[width=5cm]{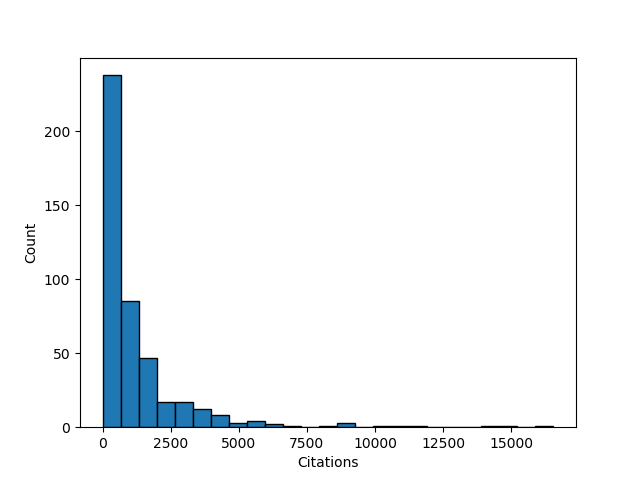}
	\includegraphics[width=5cm]{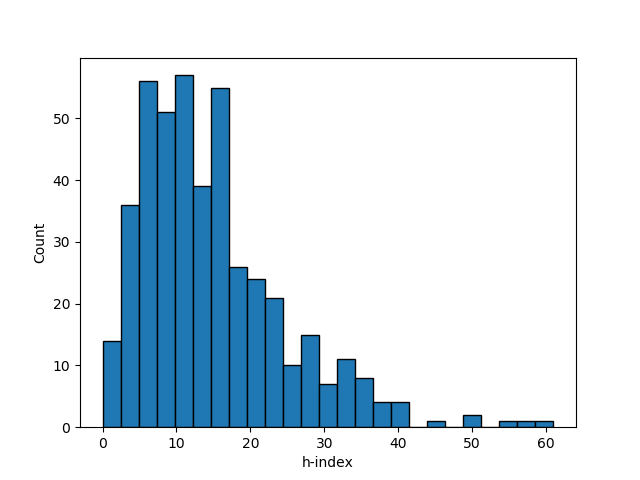}
	\includegraphics[width=5cm]{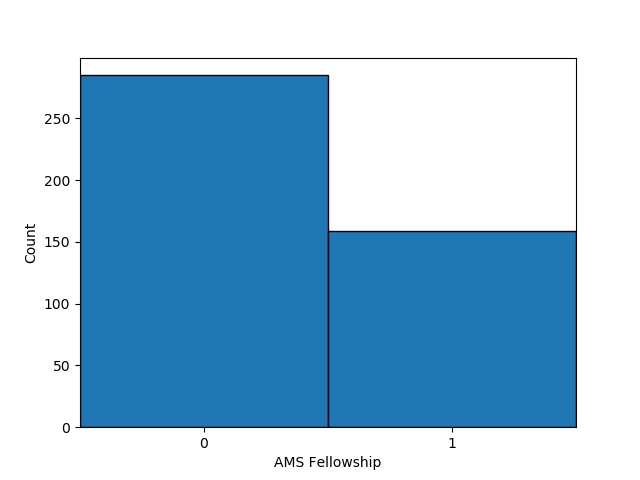}
	\includegraphics[width=5cm]{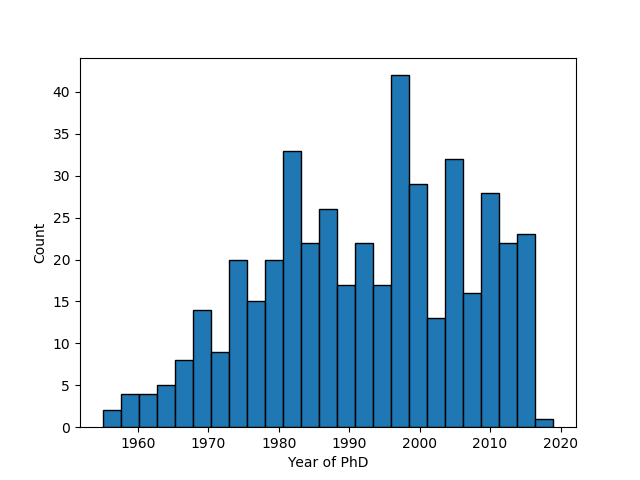}
	\caption{Histograms for each data field across all universities (n=444)}
\end{figure}

 \begin{figure}[h!]
 	\centering
 	\includegraphics[width=5cm]{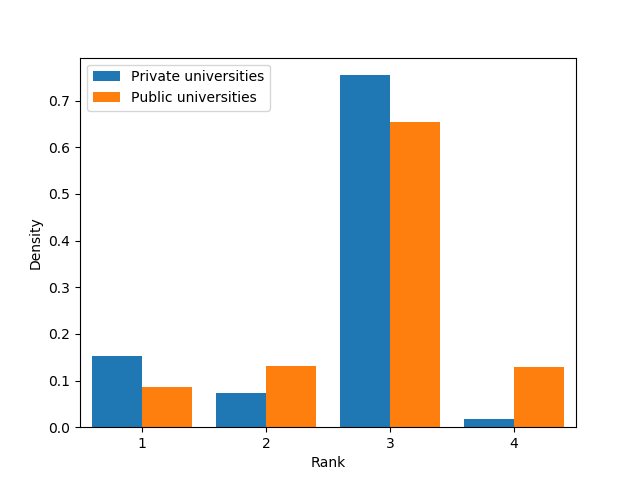}
 	\includegraphics[width=5cm]{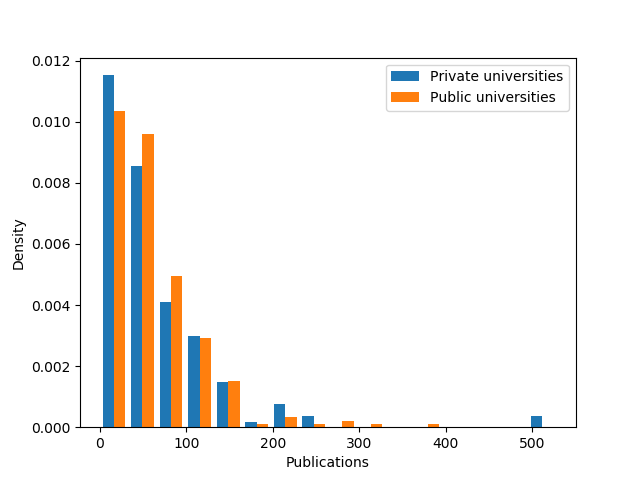}
 	\includegraphics[width=5cm]{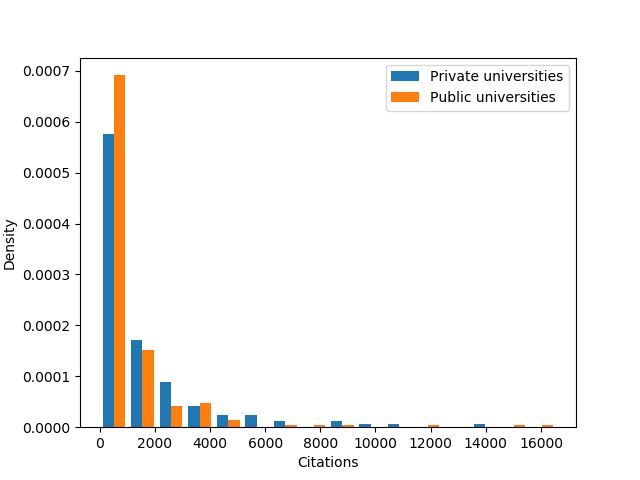}
 	\includegraphics[width=5cm]{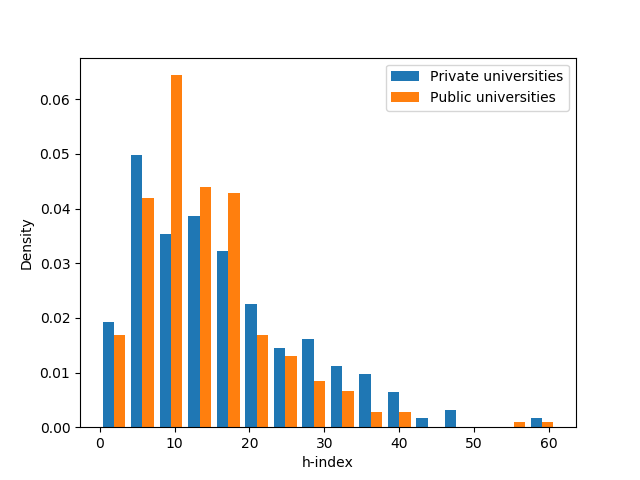}
 	\includegraphics[width=5cm]{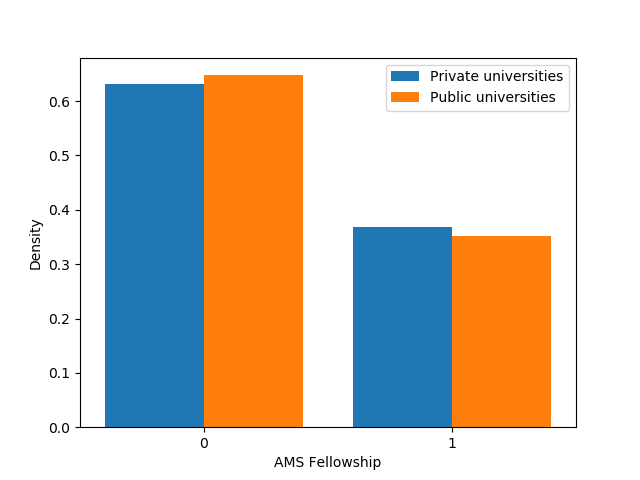}
 	\includegraphics[width=5cm]{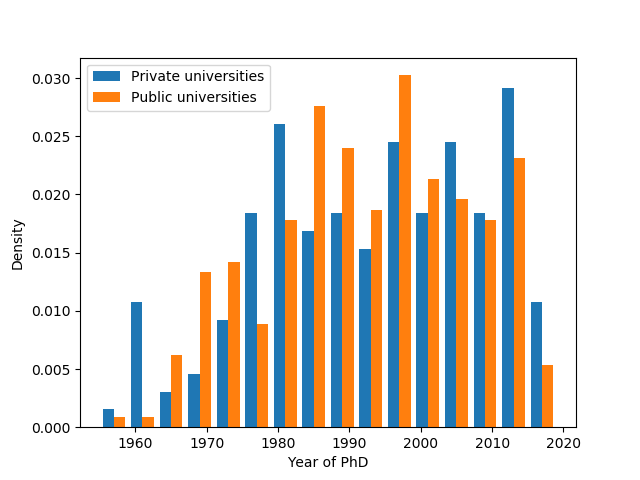}
 	\caption{Normalized histograms for each data field across all private and public universities (n=444)}
 \end{figure}

 \begin{figure}[h!]
	\centering
	\includegraphics[width=7.5cm]{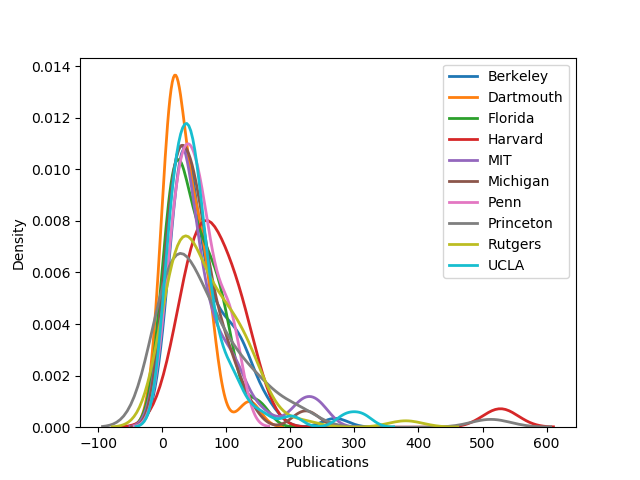}
	\includegraphics[width=7.5cm]{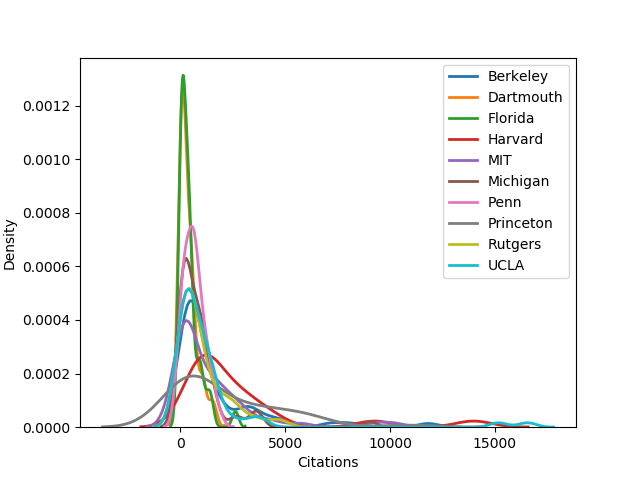}
	\includegraphics[width=7.5cm]{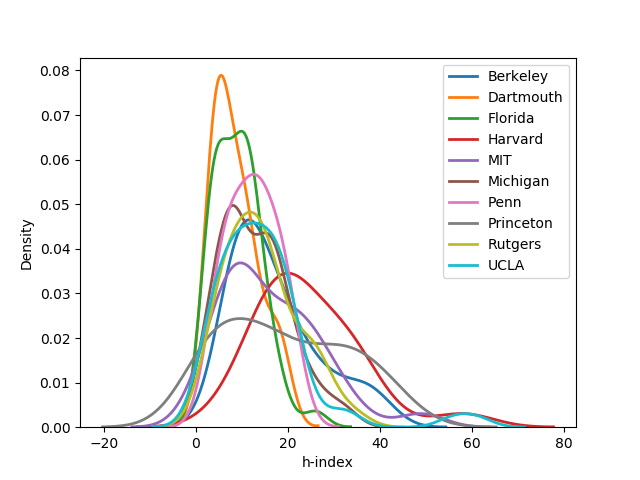}
	\includegraphics[width=7.5cm]{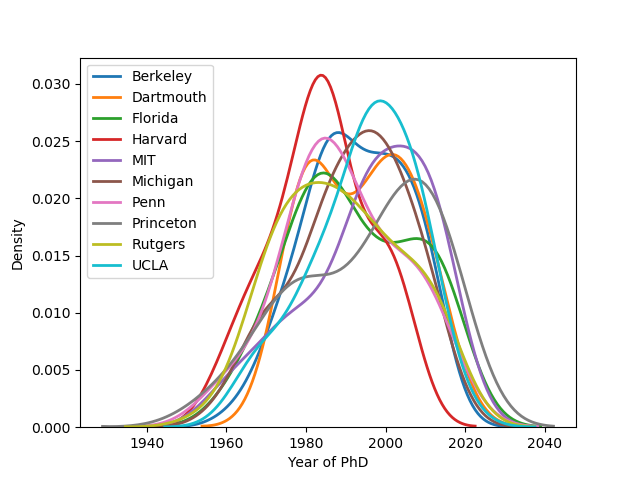}
	\caption{Kernel density estimates for each data field across each university (n=444)}
\end{figure}

We further analyzed the data for living Fields medalists as well as all Abel Prize recipients.

\begin{table}[h!]
	\centering
	\begin{tabular}{|r|c | c|}
		
		\hline	
		Field & Mean & Standard Deviation \\
		\hline
		
		Number of Publications & 130.300 & 110.607\\
		Number of Citations & 5201.800 & 5663.499\\
		\textit{h}-index & 31.050 & 16.399\\
		Year of PhD & 1982.175 & 16.522\\
		\hline	
	\end{tabular}
	\caption{Means and standard deviations across Fields medalists (n=40)}
\end{table}

\begin{table}[h!]
	\centering
	\begin{tabular}{|r|c | c|}
		
		\hline	
		Field & Mean & Standard Deviation \\
		\hline
		
		Number of Publications & 142.300 & 75.385\\
		Number of Citations & 6687.050 & 4538.117\\
		\textit{h}-index & 34.600 & 14.095\\
		Year of PhD & 1958.450 & 8.988\\
		\hline	
	\end{tabular}
	\caption{Means and standard deviations across Abel Prize recipients (n=20)}
\end{table}

\begin{table}[h!]
	\centering
	\begin{tabular}{ |c|c|c|c|c|}
		\hline	
		& Publications & Citations & \textit{h}-index & Year of PhD   \\
		\hline
		Publications & 1.000 & 0.777 & 0.807 & -0.397 \\ \hline
		Citations & 0.777 & 1.000 & 0.946 & -0.413 \\ \hline
		\textit{h}-index & 0.807 & 0.946 & 1.000 & -0.491 \\ \hline
		Year of PhD & -0.397 & -0.413 & -0.491 & 1.000 \\ \hline
	\end{tabular}
	\caption{The correlation matrix for Fields medalists (n=40)}
\end{table}

\begin{table}[h!]
	\centering
	\begin{tabular}{ |c|c|c|c|c|}
		\hline	
		& Publications & Citations & \textit{h}-index & Year of PhD   \\
		\hline
	Publications & 1.000 & 0.573 & 0.705 & -0.177 \\ \hline
	Citations & 0.573 & 1.000 & 0.944 & -0.213 \\ \hline
	\textit{h}-index & 0.705 & 0.944 & 1.000 & -0.193 \\ \hline
	Year of PhD & -0.177 & -0.213 & -0.193 & 1.000 \\ \hline
	\end{tabular}
	\caption{The correlation matrix for Abel Prize recipients (n=20)}
\end{table}

 \begin{figure}[h!]
	\centering
	\includegraphics[width=7.5cm]{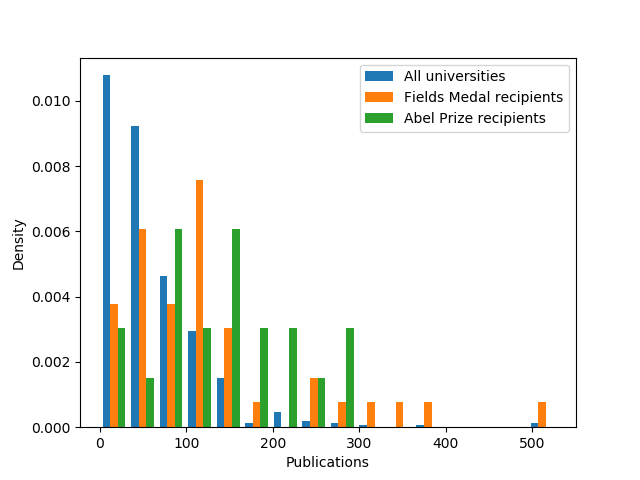}
	\includegraphics[width=7.5cm]{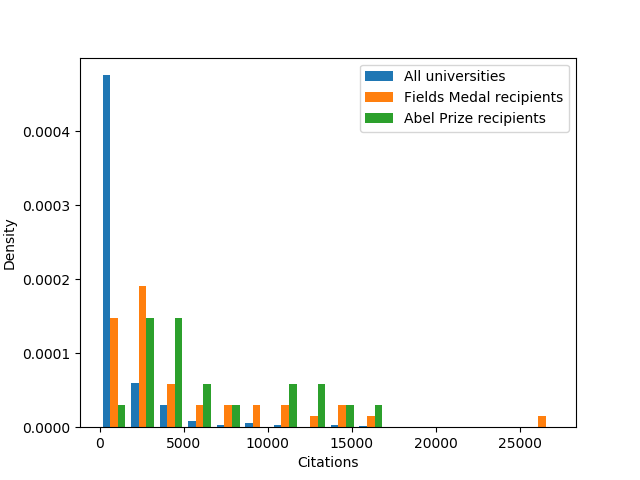}
	\includegraphics[width=7.5cm]{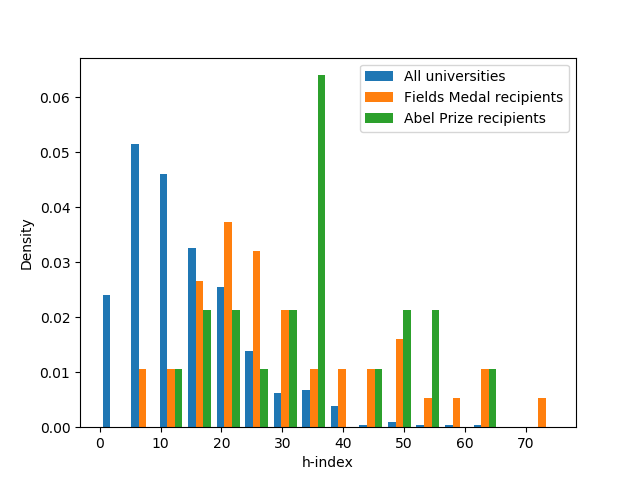}
	\includegraphics[width=7.5cm]{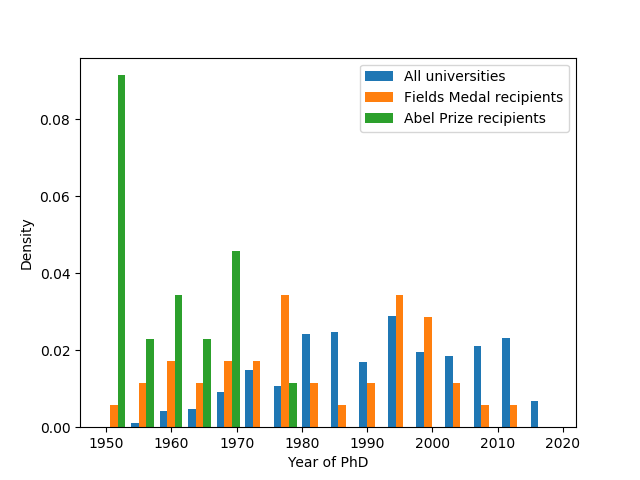}
	\caption{Normalized histograms for each data field comparing university professors (n=444) to Fields medalists (n=40) and Abel Prize reciepients (n=20)}
\end{figure}

	\section{Regression Method to Predict Rank and AMS-fellowship}
	In this section a regression method is used to attempt to predict the rank and AMS-fellow status from Number of publications, number of citations, h-index and year of PhD. However, we cannot apply these regression method directly, since the results are discrete. Specifically, our ranking goal can be regarded as a classification problem. Regression methods are often highly successful with binary classification problems. However, when we come to a multi-states classification problem, such as the four ranks we wish to predict, it seems less obvious to use these methods (except perhaps the logistic regression method). Therefore, we will make some changes on these regression method to make them can be applied on multi-classification problem. Section \ref{classificationmethod} details the regression method that we use. The rest of the section applies this methodology to predict rank and AMS-fellowship respectively of candidates. Furthermore, we will find the best combination of the four predictors, Number of publications, number of citations, h-index and year of PhD.
	
	\begin{table}[ht]
		\centering
		\label{table: table_for_regression_method}
		\begin{tabular*}{6cm}{lc}
			\hline
			 &Name of Regression Method \\
			\hline
			1 &Linear Regression\\
			2 &Logistic Regression\\
			3 &Polynomial Regression\\
			4 &RidgeCV Regression\\
			5 &Lasso Regression\\
			6 &ElasticNet Regression\\
			7 &Bayesian Ridge Regression\\
			\hline
		\end{tabular*}
	\caption{Regression Method}
	\end{table}
 
	\subsection{Classification Method}\label{classificationmethod}
		Among these method in Table 22, it is easily to use the second method, the logistic regression to deal with classification problem, even for the the problem with more than two categories. For the rest six regression methods, we simply classify the result by metrics: For each Regression methods, we can denote it by a function $F$ which maps the predictors, for example $(p_1,p_2,p_3)$ into a predicted result $py$. Since $py$ might not be the value corresponding to each category, we define a classification operator $C$ which maps each $py$ into the category whose index is the nearest value to $[py]$. Here $[x]$ represents the floor function; ie, the largest integer less than $x$.  Therefore, the classification method can be represented by 
		\[ C\circ F \].
		
	\subsection{Data sets}
	In this section, we randomly choose 70 percents of the whole data set to be trained and the rest to be tested.
	\subsection{Evaluate Index}
	In order to evaluate the fitness of a method, we use two different index. One is the average degree of deviation, $ADD$, and the other is the accuracy rate, $AR$. For the predicted result of test data $\{py_i\}_{i=1}^N$ and the real test data $\{y_i\}_{i=1}^N$, we define the variance $ADD$ and $AR$ to be:
   \begin{equation}
	ADD:=\frac{\sum_{i=1}^{N} \left|py_i - y_i\right|}{N}
\end{equation}
\begin{equation}
	AR:=\frac{\sum_{i=1}^N \delta(py_i,y_i)}{N}
	\end{equation}
	
	where the $\delta(py,y)$ is defined as following
\begin{equation}
	\delta(py,y):=\left\{
	\begin{array}{ll}
	1 &if\ py = y\\
	0 &others
	\end{array}
	\right.
	\end{equation}
A higher value of $AR$ corresponds with a higher accuracy of the predicted result, and a lower value of $ADD$ with a more robust regression method.
	
	\subsection{Prediction Result}
	In this section, we examine the results of different regression methods and different combination of these four predictors, corresponding to different university. Since there are seven different regression methods, fifteen combinations of predictors, ten schools, two features to be predicted and two evaluation factors, we do not list the entirety of the calculated results in this paper. Instead, for each school, we report the regression method and predictors combination with highest $AR$ or lowest $ADD$.  
	
	We denote different combination of predictors by different index as following:
	\begin{table}[ht]
		\centering
		\label{table: table_for_index_to_combination}
		\begin{tabular*}{12cm}{llll}
			\hline
			Index &Combination of predictors &Index &Combination of predictors\\
			\hline
			1 &Number of publications &8 &2 and 3\\
			2 &Number of citations &9 &2 and 4\\
			3 &h-index &10 &3 and 4\\
			4 &year of PhD &11 &1 and 2 and 3\\
			5 &1 and 2 &12 &1 and 2 and 4\\
			6 &1 and 3 &13 &1 and 3 and 4\\
			7 &1 and 4 &14 &1 and 2 and 3 and 4\\
			\hline
		\end{tabular*}
		\caption{Index of different combination}
	\end{table}

We introduce our prediction method with the Berkeley data as an example, then the rest of the section will list the results.  
	
	\subsubsection{Berkeley}
	As shown in Tables 24, 25, 26 and 27, we can conclude that:
	\begin{itemize}
		\item To predict the rank for $Berkeley$, we can use $LnR$, $RR$, $LR$, $ENR$ and $ByR$ method and the predictor combination can be 1 to 14 except from 11 and 14. However, the low $AR$ value shows that it doesn't seem a good method to use regression to predict the rank for Berkeley. 
		\item To predict the AMS-fellowship status for $Berkeley$, we use the method and predictor pair $(PoR, 5)$ which has both the highest $AR$ value and the lowest $ADD$ value. (It seems a coincident.) Moreover, we can write down the formula for this $(PoR, 5)$ pair. 
		\[
		f(Np,Nc\footnote{Np is the abbreviation of Number of Publications and Nc is the abbreviation of Number of citations})=(81-179Np-143Np^2+31Np^3)\cdot(-678-3730Nc+625Nc^2-1846Nc^3)
		\]
		With $f(Np,Nc)$, we can predict $AMS-fellowship$ by
	\begin{equation}
		AMS = \left\{
		\begin{array}{ll}
		1 &f(Np, Nc) \ge 1\\
		0 &f(Np, Nc) <1\\
		\end{array}
		\right.
	\end{equation}
	\end{itemize}
	\begin{table}[ht]
		\centering
		\label{table: table_for_ADD_rank}
		\begin{tabular*}{9cm}{lccccccc}
			\hline
			&LnR &LgR &PoR &RR &LR &ENR &ByR\\
			\hline
			1 &0.89 &1.22 &1.33 &0.89 &0.89 &0.89 &0.89 \\
			2 &0.89 &1.22 &1.33 &0.89 &0.89 &0.89 &0.89 \\
			3 &0.89 &1.22 &1.33 &0.89 &0.89 &0.89 &0.89 \\
			4 &0.89 &1.33 &1.33 &0.89 &0.89 &0.89 &0.89 \\
			5 &0.89 &1.22 &1.33 &0.89 &0.89 &0.89 &0.89 \\
			6 &0.89 &1.22 &1.33 &0.89 &0.89 &0.89 &0.89 \\
			7 &0.89 &1.22 &1.33 &0.89 &0.89 &0.89 &0.89 \\
			8 &0.89 &1.22 &1.33 &0.89 &0.89 &0.89 &0.89 \\
			9 &0.89 &1.22 &1.22 &0.89 &0.89 &0.89 &0.89 \\
			10 &0.89 &1.00 &1.33 &0.89 &0.89 &0.89 &0.89 \\
			11 &0.94 &1.22 &1.33 &0.94 &0.94 &0.89 &0.89 \\
			12 &0.89 &1.22 &1.22 &0.89 &0.89 &0.89 &0.89 \\
			13 &0.89 &1.22 &1.22 &0.89 &0.89 &0.89 &0.89 \\
			14 &0.94 &1.22 &1.22 &0.89 &0.94 &0.89 &0.89 \\
			\hline
		\end{tabular*}
		\caption{ADD of the prediction results for rank, Berkeley}
	\end{table}

	\begin{table}[ht]
		\centering
		\label{table: table_for_ADD_AMS}
		\begin{tabular*}{9cm}{lccccccc}
			\hline
			&LnR &LgR &PoR &RR &LR &ENR &ByR\\
			\hline
			1 &0.61 &0.61 &0.61 &0.61 &0.61 &0.61 &0.61 \\
			2 &0.56 &0.44 &0.39 &0.56 &0.56 &0.56 &0.56 \\
			3 &0.61 &0.44 &0.61 &0.61 &0.61 &0.61 &0.61 \\
			4 &0.61 &0.61 &0.61 &0.61 &0.61 &0.61 &0.61 \\
			5 &0.50 &0.50 &0.28 &0.50 &0.50 &0.50 &0.56 \\
			6 &0.61 &0.39 &0.61 &0.61 &0.61 &0.61 &0.61 \\
			7 &0.61 &0.61 &0.61 &0.61 &0.61 &0.61 &0.61 \\
			8 &0.56 &0.44 &0.61 &0.56 &0.56 &0.56 &0.56 \\
			9 &0.56 &0.44 &0.50 &0.56 &0.56 &0.56 &0.56 \\
			10 &0.61 &0.44 &0.61 &0.61 &0.61 &0.61 &0.61 \\
			11 &0.50 &0.44 &0.61 &0.50 &0.50 &0.50 &0.56 \\
			12 &0.50 &0.50 &0.50 &0.50 &0.50 &0.50 &0.56 \\
			13 &0.56 &0.44 &0.50 &0.56 &0.56 &0.56 &0.56 \\
			14 &0.50 &0.44 &0.50 &0.50 &0.50 &0.50 &0.56 \\
			\hline
		\end{tabular*}
		\caption{ADD of the prediction results for AMS-fellowship status, Berkeley}
	\end{table}

	\begin{table}[ht]
		\centering
		\label{table: table_for_AR_rank}
		\begin{tabular*}{9cm}{lccccccc}
			\hline
			&LnR &LgR &PoR &RR &LR &ENR &ByR\\
			\hline
			1 &0.11 &0.06 &0.00 &0.11 &0.11 &0.11 &0.11 \\
			2 &0.11 &0.06 &0.00 &0.11 &0.11 &0.11 &0.11 \\
			3 &0.11 &0.06 &0.00 &0.11 &0.11 &0.11 &0.11 \\
			4 &0.11 &0.00 &0.00 &0.11 &0.11 &0.11 &0.11 \\
			5 &0.11 &0.06 &0.00 &0.11 &0.11 &0.11 &0.11 \\
			6 &0.11 &0.06 &0.00 &0.11 &0.11 &0.11 &0.11 \\
			7 &0.11 &0.06 &0.00 &0.11 &0.11 &0.11 &0.11 \\
			8 &0.11 &0.06 &0.00 &0.11 &0.11 &0.11 &0.11 \\
			9 &0.11 &0.06 &0.06 &0.11 &0.11 &0.11 &0.11 \\
			10 &0.11 &0.11 &0.00 &0.11 &0.11 &0.11 &0.11 \\
			11 &0.11 &0.06 &0.00 &0.11 &0.11 &0.11 &0.11 \\
			12 &0.11 &0.06 &0.06 &0.11 &0.11 &0.11 &0.11 \\
			13 &0.11 &0.06 &0.06 &0.11 &0.11 &0.11 &0.11 \\
			14 &0.11 &0.06 &0.06 &0.11 &0.11 &0.11 &0.11 \\
			\hline
		\end{tabular*}
		\caption{AR of the prediction results for rank, Berkeley}
	\end{table}
	
	\begin{table}[ht]
		\centering
		\label{table: table_for_AR_AMS}
		\begin{tabular*}{9cm}{lccccccc}
			\hline
			&LnR &LgR &PoR &RR &LR &ENR &ByR\\
			\hline
			1 &0.39 &0.39 &0.39 &0.39 &0.39 &0.39 &0.39 \\
			2 &0.44 &0.56 &0.61 &0.44 &0.44 &0.44 &0.44 \\
			3 &0.39 &0.56 &0.39 &0.39 &0.39 &0.39 &0.39 \\
			4 &0.39 &0.39 &0.39 &0.39 &0.39 &0.39 &0.39 \\
			5 &0.50 &0.50 &0.72 &0.50 &0.50 &0.50 &0.44 \\
			6 &0.39 &0.61 &0.39 &0.39 &0.39 &0.39 &0.39 \\
			7 &0.39 &0.39 &0.39 &0.39 &0.39 &0.39 &0.39 \\
			8 &0.44 &0.56 &0.39 &0.44 &0.44 &0.44 &0.44 \\
			9 &0.44 &0.56 &0.50 &0.44 &0.44 &0.44 &0.44 \\
			10 &0.39 &0.56 &0.39 &0.39 &0.39 &0.39 &0.39 \\
			11 &0.50 &0.56 &0.39 &0.50 &0.50 &0.50 &0.44 \\
			12 &0.50 &0.50 &0.50 &0.50 &0.50 &0.50 &0.44 \\
			13 &0.44 &0.56 &0.50 &0.44 &0.44 &0.44 &0.44 \\
			14 &0.50 &0.56 &0.50 &0.50 &0.50 &0.50 &0.44 \\
			\hline
		\end{tabular*}
		\caption{AR of the prediction results for AMS-fellowship status, Berkeley}
	\end{table}
	
	With this method, we can write down the best method-predictors pair, the corresponding AR and ADD value to it and the prediction formula for each school.
	\subsubsection{Dartmouth}
	In Dartmouth, all the regression method and predictor combinations have a 100\% accuracy rating for predicting AMS Fellowship status. However, these methods were worse for the rank prediction. The best prediction pair has only 29\% AR and 0.71 AD.  
	\subsubsection{Florida}
	In Florida, there is also 100\% accuracy in AMS status. The best prediction pair for professorship rank has 31\% AD and 0.69 AD.
	\subsubsection{Harvard}
	In Harvard, the best prediction pairs for AMS status are $(8,LgR)$, $(13,LgR)$, $(9,PoR)$, $(12,PoR)$, $(13,PoR)$ and $(14, PoR)$ with 67\% AR and 0.5 ADD.
	\subsubsection{Michigan}
	In Michigan, the best prediction pairs for AMS status are $(3,LgR)$, $(8,LgR)$, $(10,LgR)$, $(13, LgR)$ and $(12,RR)$ with 79\% AR and 0.21 ADD. The best prediction pair for rank is $(9,RR)$ with 32\% AR and 0.68 ADD.
	\subsubsection{MIT}
	In MIT, the best prediction pairs for AMS status are $(1, LgR)$, $(3,LgR)$, $(7,LgR)$ and $(10,LgR)$ with 69\% AR and 0.31 ADD. The best prediction pairs for rank are $(7,LgR)$, $(11,LgR)$ and $(14,LgR)$ with 12\% AR and 0.88 ADD. 
	\subsubsection{Upenn}
	In Upenn, the best prediction pair for AMS status is $(1,PoR)$ with 62\% AR and 0.38 ADD. The best prediction pairs for rank are $(1,LgR)$, $(2,LgR)$, $(5,LgR)$, $(6,LgR)$, $(7,LgR)$, $(8,LgR)$, $(9,LgR)$, $(11,LgR)$, $(12,LgR)$, $(13,LgR)$, $(14,LgR)$, $(9, PoR)$, $(12, PoR)$, $(13, PoR)$ and $(14, PoR)$ with 12\% AR and 0.88 ADD.
	\subsubsection{Princeton}
	In Princeton, the best prediction pair for AMS status is $(4,LgR)$ with 92\% AR and 0.08 ADD. The best prediction pair for rank has 23\% AR and 0.85 ADD. 
	\subsubsection{Rutgers}
	In Rutgers, the best prediction pairs for AMS status are $(3,LgR)$, $(11,LgR)$, $(13,LgR)$ and $(14,LgR)$ with 89\% AR and 0.11 ADD. The best prediction pairs for rank are $(14,LnR)$, $(14, RR)$ and $(14,LR)$ with 78\% AR and 0.22 ADD.
	\subsubsection{UCLA}
	In UCLA,the best prediction pairs for AMS status are $(5,PoR)$ and $(11,PoR)$ with 61\% AR and 0.39 ADD. All the prediction pairs for rank work far worse with only 6\% AR and 0.94 ADD at most.

\section{Artificial Neural Network Model}
Artificial neural network (ANN), or deep learning, is a specific subfield of machine learning and a new method on learning representations from data which puts an emphasis on learning successive ``layers'' of increasingly meaningful representations. The name ``neural network'' is from brain science, however, ANN is merely a mathematical framework for learning representations from data. A deep network can be imagined as a multi-stage information distillation operation, where information goes through successive filters and comes out increasingly ``purified'', i.e., useful with regard to some task. There are rich literatures on ANN, and more broadly, machine learning. We refer interested readers to \cite{C} and \cite{G} for more theoretical backgrounds and hands-on skills on these topics. 

As Francous Chollet says, machine learning is an art rather than a science. There are no definite rules telling one what choices of architectures, hyperparameters, etc. will lead to the optimal results. Hence, we would like to explore Artificial Neural Network (ANN) models with different settings in this section. 

For the study on ANN models, we mainly use the {\tt keras} module in python. This is a high level API utilizing {\tt tensorflow} as its backend. There are two kinds of models in {\tt keras}, sequential model and function API. The first one is more popular and satisfies most needs. Functional API can help one construct any network, i.e. a graph where each node of the graph is a layer in the model. Each layer consists of a few hidden units, or neurons, in either model. There are several options to connect adjacent layers, the most popular one being {\tt Dense}, i.e., a unit in a layer is connected to all units in its adjacent layer(s). Other types of layers include locally-connected layers, recurrent layers, convolutional layers, embedding layers, merge layers, normalization layers and noise layers, etc.

\subsection{Prediction on the Rank with ANN}
Artificial neural networks (ANNs) are used most often to extract complex 
relationships within a data set. Considering different departments usually have different standards for promotion and various rank structures, i.e., no distinguished professor at Princeton and Harvard as mentioned, we try the prediction with Rutgers data set as an example. The studies on other departments are similar and are left to interested readers.

While our data set is currently small,  it is still interesting to explore how well an ANN  classifies our data. Even with such a small data set, one can see a decent prediction power. More precisely, in this section, we will work with the following problem: given a professor, who is represented by a list of length three  containing the number of publications, the number of  citations on these publications, and h-index, we would like to predict what rank this professor has. In order to use an ANN for this task, our final model should output a vector which has length of the number of possible rankings whose elements are between zero and one and whose entries sum to 1. We will begin with a simple linear classifier, and  after experimenting with this simple model, but future research will address if adding non-linearity through a second layer can increase the accuracy. 

One of the simplest forms of an ANN is a one-layer linear classifier. 
We shall follow the article {\tt http://cs231n.github.io/neural-networks-case-study/\#linear}
from the Stanford cs231n course with several modifications. The model
explained in this article is known as a \textit{soft-max linear classifier}.
In this model, our data undergoes a linear transformation from some 
$k$-dimensional real space to $N$-dimensional real space, where $k$
is the number of descriptive features of the data and $N$ is the 
number of target features. We interpret this $N$-dimensional vector
as a list of unnormalized log probabilities, and we apply the soft-max 
function which element-wise exponentiates and normalizes this vector to obtain a list
of probabilities. 

We will train our neural network with hand-labeled (by the Rutgers Mathematics 
promotion committee) data, which is a list of professors and current rankings
in the format [descriptive feature 1, descriptive feature 2,\ldots, descriptive feature k, rank].
Descriptive features may be chosen from the following: number of citations,
number of publications, $h$-index, AMS status, and year of receiving PhD.
Before training the ANN, we do the following preprocessing on our training set data:
For each descriptive feature $F$, we transform $F$ so that it has mean 
zero and standard deviation one. In addition, since our model predicts 
probabilities, we convert the number professor rank into a length-four
vector (probability distribution), which is a one at position $i$ if 
the professor is of rank $i$ and zero otherwise. This is known as a 
\textit{one-hot} encoding of the target feature. Using this 
encoding of the target feature, we can compute how far wrong our model's
current prediction is from the truth. To this end, we use the cross-entropy loss function.
For two probability distributions $p$, the true distribution, and $q$, the test distribution, on a base set
 $X$, the
cross-entropy $L(p,q)$ is defined as 
$$L(p,q)=\sum_{x\in X} -p(x)\log(q(x)).$$
Using our one-hot encoding of the target feature, our loss for a single
piece of data is thus
$$-\log(q(x_i)),$$
where $i$ is the correct label for this piece of data. We sum over all
of the training data to get the loss for a single iteration (epoch) of training.
Using the loss function, we back-propagate the error after each epoch
 to update the weights of our ANN. In addition, we also update the ANN weights with 
 a small amount of \textit{regularization}, which keeps the weights closer to zero.
 The purpose of this is to prevent over-fitting our data and to encourage use of 
 all target features by the ANN.
 
We start with training the neural net using all available numerical descriptive features other
than salary. We permute the data after extracting the relevant fields and take the first 45 
entries to train the ANN. The rest we set aside for testing. After experimenting with hyper-parameters, 
we find that the network seems to converge after 200 epochs. Below is a plot of the 
cross-entropy loss for each epoch of training on this data set.

\begin{figure}[ht]
\centering
\includegraphics[scale=0.6]{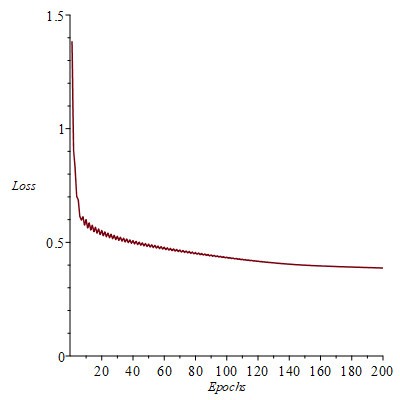}
\caption{The cross-entropy loss for each epoch of training}
\end{figure}

To test the trained network, we let the network's prediction of a given professor rank be the argmax
of the list of probabilities. We may now evaluate the accuracy on the training set and find that the ANN predicts professor rank 
correctly 12 out of 14 times! The list of predictions by the network is [4, 4, 4, 4, 2, 4, 3, 4, 3, 3, 4, 3, 1, 3], and the 
correct rankings are [4, 4, 4, 4, 2, 4, 4, 4, 3, 3, 4, 3, 1, 2]. 

It is interesting to see how our model performs with using fewer descriptive features. It seems that the number of 
publications, the number of citations, and the h-index are particularly important criteria, so we use these to train the ANN. Using the same hyper-parameters,
the neural net trains well after 200 epochs. The loss curve is similar, yet we find that 
the ANN predicts the ranking correctly only 9 out of 14 times. The list of predictions is [3, 2, 3, 1, 3, 3, 4, 3, 3, 4, 4, 3, 4, 1]
in comparison to the true rankings [3, 3, 3, 1, 2, 3, 4, 4, 4, 4, 4, 3, 4, 2]. It is instructive to plot the data to see how it is 
clustered. Below is a plot of the three-dimensional data where each points color represents the ranking of the professor: magenta corresponds to assistant professor, blue to associate professor, green to professor, and orange to distinguished professor.

\begin{figure}[ht]
\centering
\includegraphics[scale=0.7]{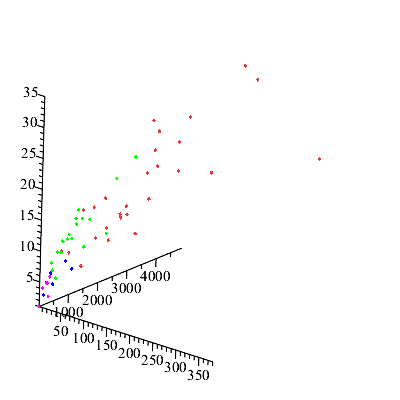}
\caption{The three-dimensional data where each points color corresponds to the ranking of the professor.}
\end{figure}

One can see that the data is not linearly separable, and in fact does not seem to be separable by any simple non-linear model.
Future investigation could include finding a small number of parameters which allow for a linear separation of the data 
or seeing how well a non-linear model predicts professor rankings.

\subsection{Exploration on Math Faculty Data with ANN}
To use the entire data set, we attempt to predict AMS status, as this is a feature that is common across all departments and not subjective to internal departmental policy. Hence, we can use all the data of publication, citation and $h$-index information to predict whether a professor is an AMS fellow.

Since there are already numerous literatures on how to tune an ANN model and how to find the ``best" hyperparameters, we merely give an example of code here. Following is an example of the algorithm. We use python's {\tt keras} package to explore the prediction. 

\begin{figure}[ht]
\centering
\includegraphics[scale=0.64]{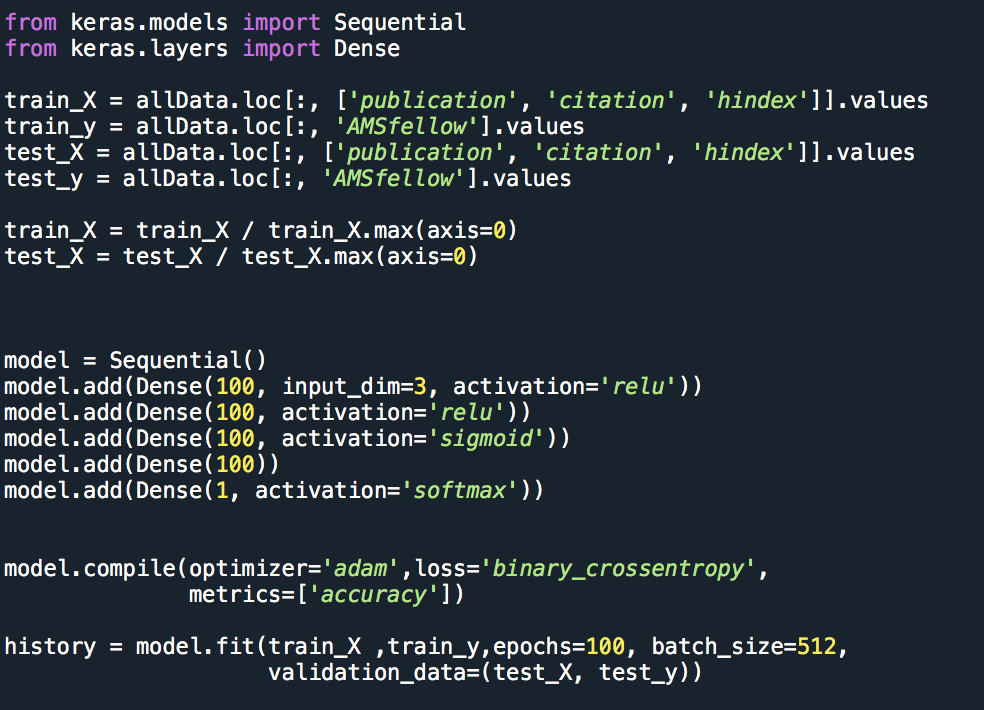}
\caption{Keras Package for ANN}
\end{figure}

By preprocessing the input data, adding regularization, trying different architecture and activation functions, doing a grid search for hyperparameters and choosing suitable metrics to evaluate the model, it would be very promising to have a great precision in the prediction. The tuning and refining process is left to interested readers.

\section{Unsupervised Clustering For Predictive Analysis}

The method used in this section is an unsupervised clustering algorithm developed by the 2015 UCLA Applied Math REU Hyperspectral Imagery research team \cite{UCLA2015},\cite{UCLA2015Paper}. It was chosen because it was designed specifically to sort large sets of data into a relatively small number of sorted groups, with no prior information or training data needed.  The following terminology will be borrowed from the hyperspectral lexicon: each sorted group is called a \textit{cluster}, and the average vector of a cluster is its \textit{centroid}. 

In the context of hyperspectral imagery, the NLTV algorithm is notable because there are very few robust unsupervised algorithms. Here, it is the lack of necessity of training data which makes it an interesting clustering method to apply: while data can be collected from universities across the country, there is no guaranteed standard of departmental promotions, which means each university ought to be treated separately, and as such that does not provide much training data for a neural network. 

\subsection{The Algorithm}

The core of the sorting algorithm comes from the minimization of an energy functional \begin{equation}
E(u)=\parallel\bigtriangledown u \parallel_{L1} + \lambda \langle u, f \rangle,
\end{equation}
where $u:\Omega\rightarrow[0,1]^n$ is the \textit{labeling function} on the data, $n$ is the number of clusters it is being sorted into, and $\Omega$ is the domain of the data, and \textit{f} is a fidelity function. The inspiration comes from the imaging process technique of total variation introduced by Rudin et al in 1992 \cite{totalvariation1} for noise reduction, which corresponds to the minimization of the gradient of $u$. In highly noisy images or datasets where adjacent pixels do not matter, simply calculating the gradient directly does not give as pertinent information. Therefore, we turn to the theory of nonlocal operators introduced by \cite{zhou1},\cite{zhou2}, Zhou and Sch\"{o}lkopf and adapted to image processing by Osher and Gilboa \cite{oshernonlocal}.

Let $\Omega$ be a region in $\mathbb{R}^k$, and $u : \Omega \to \mathbb{R}$ be a real function. Then the non-local derivative is defined as \begin{equation} \frac{\partial u}{\partial y}(x):=\frac{u(y)-u(x)}{d(x,y)},\hspace{10pt} \text{for all } x, y \in \Omega \end{equation} where $d$ is a positive distance between $x$ and $y$. With the following non-local weight defined as \ref{weight}, we can re-write the non-local derivative as \ref{nld}. \begin{equation} w(x,y)=d^{-2}(x,y) \label{weight}\end{equation} \begin{equation} \frac{\partial u}{\partial y}(x) = \sqrt{w(x,y)}(u(y)-u(x)) \label{nld}\end{equation} 
Then the non-local gradient $\bigtriangledown_{w} u$ for $u \in L^{2}(\Omega)$ as a function from $\Omega$ to $L^{2}(\Omega)$ is the collection of all partial derivatives \begin{equation} \bigtriangledown_w u(x)(y)=\frac{\partial u}{\partial y}(x) = \sqrt{w(x,y)}(u(y)-u(x)). \end{equation} 

Note that here, ``distance" can either refer to the standard Euclidean distance \begin{equation}
d(x,y)=\sqrt{\sum_{i=1}^k(x_i-y_i)^2},
\end{equation}
the cosine distance \begin{equation}
d(x,y)=1-\frac{x\cdot y}{||x||||y||},
\end{equation}
or a linear combination of them.

The non-local energy functional we are trying to minimize takes the form of \begin{equation}
E(u)=\parallel\bigtriangledown_w u \parallel_{L1} + \lambda\sum_{i=1}^n|u_i(x)g(x)-c_i|^2\label{EF},
\end{equation}
where $\parallel\bigtriangledown_w u \parallel_{L1}$ is the $L^1$ norm on the space $L^2(\Omega,L^2(\Omega))$ defined as
\begin{equation}
	\parallel v \parallel_{L^1}:=\int_\Omega \parallel v(x) \parallel_{L^2} dx=\int_\Omega\Big|\int_\Omega v(x)(y)^2dy\Big|^{\frac{1}
		{2}}dx
\end{equation}
and the fidelity function is explicitly given by $\lambda\sum_{i=1}^n|u_i(x)g(x)-c_i|^2$, where $g(x)$ is the datapoint and $c_i$ is the $i$th cluster centroid. We explicitly discretize the labeling function and nonlocal operators, $u=(u_1,u_2,\ldots,u_n)$ is a matrix of size $m\times n$, where $m$ is the number of datapoints and $n$ is the number of clusters. Each $u_i$ takes values between 0 and 1, $\sum_{i=1}^n u_{ki}=1$ for all $k\in{1,...,m}$. Then $(\bigtriangledown_wu_l)_{i,j}=\sqrt{w_{i,j}}((u_l)_j-(u_l)_i)$ is the nonlocal gradient of $u_l$; $(\text{div}_wv)_i=\sum_j\sqrt{w_{i,j}}v_{i,j}-\sqrt{w_{j,i}}v_{j,i}$ is the divergence of $v$ at $i$-th datapoint; and the discrete $L^1$ norm of $\bigtriangledown_wu_l$ are defined as: \begin{equation}\parallel \bigtriangledown_wu_l \parallel_{L^1}=\sum_i\left(\sum_j\left(\bigtriangledown_wu_l\right)_{i,j}^2\right)^{\frac{1}{2}}.\end{equation}

The functional \ref{EF} is convex, so a global minimum exists. However, calculating $\parallel\bigtriangledown u \parallel_{L^1}$ via gradient descent involves calculating div($\frac{\bigtriangledown u}{\mid \bigtriangledown u \mid}$), which is highly unstable because $\mid \bigtriangledown u \mid$ can be equal to zero. In 2011, Chambolle and Pock introduced a first-order primal dual algorithm, which they proved converged to a saddle point with a rate of $O(1/N)$ in finite dimensions for the complete class of convex problems \cite{primaldual}. This was used as an inspiration to craft a saddle point solution with respect to $u,$$\overline{u}$, and $p$. Full motivation and description can be found in \cite{UCLA2015}, \cite{UCLA2015Paper},\cite{WeiThesis}. The algorithm is as follows:

\begin{mdframed}
	\textbf{Primal-Dual Iterations}\\
	$\bullet$ Iterations $(n> 0)$: Update $u^n, p^n, \bar{u}^n$ as follows:\\
	\begin{equation}
		\nonumber
		\begin{cases}
			p^{n+1}=\text{proj}_P(p^n+\sigma \bigtriangledown_w \bar{u}^n)\\
			u^{n+1}=\arg\min_u \delta_U(u)+\frac{1}{2}\parallel (I+\tau F)^{\frac{1}{2}}u-(I+\tau 
			F)^{-\frac{1}{2}}(u^n+\tau \text{div}_wp^{n+1})\parallel ^2\\
			\bar{u}^{n+1}=u^{n+1}+\theta (u^{n+1}-u^n)\\
		\end{cases}
	\end{equation}
\end{mdframed}
where $F$ is the discretized fidelity function matrix with the inbuilt weight $\lambda$.

The overall sorting algorithm is then:

\begin{mdframed}
\textbf{Nonlocal Total Variation Unsupervised Clustering}\\
$\bullet$ Initiate parameters.\\
$\bullet$ Calculate weight matrix.\\
$\bullet$ Set $n$ random datapoints as the first iteration of centroids, set\\ \textcolor{white}{...}$u^0=\overline{u}^0=$ Matrix(m,n,$1/m$) and $p^0$ zeroed out.\\
\textbf{while} not converge do \\
\textcolor{white}{......}\textbf{Inner Loop:} Primal Dual Algorithm to find minimizing $u$.\\
\textcolor{white}{......}\textbf{Outer Loop:} Threshold $u$ into an assignment function, and use the new sorting of the data to update the centroids.\\
\textbf{end}
\end{mdframed}

The NLTV algorithms was originally retooled for clustering of mathematical data by the authors in \cite{DAMF}, which only used data from Rutgers University for sorting but included data on salaries. There are two main changes between the algorithm written for this paper, and the algorithm developed in 2015. Firstly, the calculation of the weight matrix is done directly between all datapoints in this project; the original hyperspectral algorithm used a ``patch" distance to filter for noise, and employed an approximate nearest neighbor search to save computational time. In the original algorithm, a smart simplex clustering method instead of directly thresholding was developed for the hyperspectral with inspiration from \cite{hnmf}. The final thresholding process was not used in the analysis of the data as it makes the outerloop of the algorithm far more computationally expensive, for no increase in convergence time in a dataset as clean as this one.

There are a number of parameters involved in the algorithm, but the two most vital ones are $\lambda$, which determines the weight given to the minimization of the fidelity function vs the gradient of $u$, and the choice of Euclidean vs Cosine distance for the creation of the weight matrix and fidelity distance calculations. The value for $\lambda$ ought to be comparatively large to prioritize tight sorting. Euclidean vs Cosine vs a linear combination is something that should be tailored to the dataset, as some of the fields (ie h-index or AMS Fellow: 0/1) have a smaller range of values, and some of the fields (ie number of citations or year of PhD) have a much larger range of values, so that field does not dominate. 

\subsection{Results}
\begin{figure} 
\includegraphics[width=\textwidth]{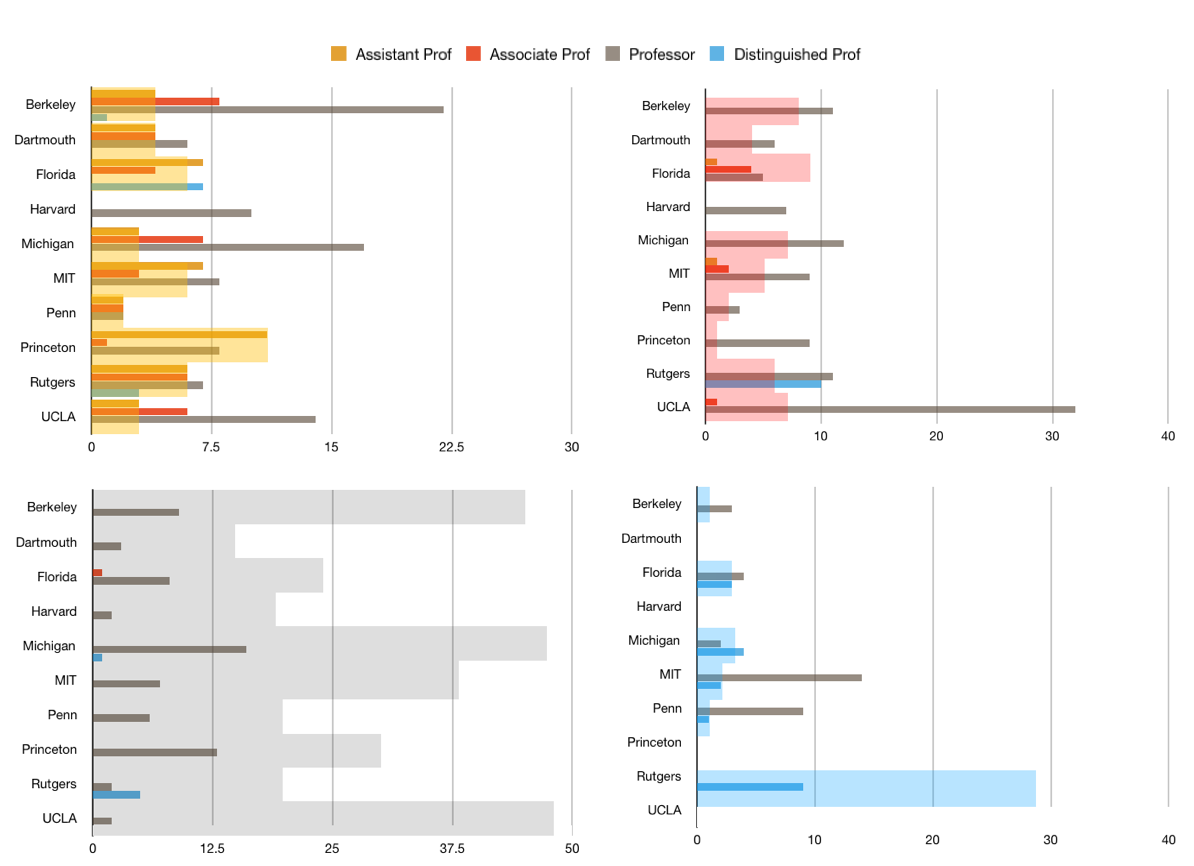}   
\caption{Sorted Clusters Vs Ground Truth}\label{cl}  
\end{figure}

The individual results for the average three or four centroids of ten universities are listed in the charts below. `Rank' indicates 1 for Assistant Professor, 2 for Associate Professor, 3 for Professor, and 4 for Distinguished Professor, which `AMS' denotes 1 for AMS Fellow, and 0 if not.  Figure \ref{cl} gives a secondary direct visual of the ``accuracy" of each cluster by denoting the actual ranks of each professor sorted into the associated centroid. Some universities did not have Distinguished Professors, and hence the data was sorted into three clusters instead of four. Harvard only had Professors, and so was sorted into three clusters. Parameters 'cosine' indicates Euclidean weight $10^{-10}$, Cosine weight 1, $\lambda=1$, and 'mixed' indicates Euclidean weight 1, Cosine weight $10^2$, $\lambda=10^4$.

The general pattern of the results is as follows: the NLTV clustering algorithm is usually able to pick out the extremes correctly (i.e. placing all of Rank 1 or Rank 4 in the same cluster); however, the extremely large variance in the Professor / Rank 3 Category means that oftentimes multiple Professor clusters would form instead of the desired ranking.

\textbf{Berkeley}, \textit{Parameters:} Mixed.

\begin{tabular}{|c ||c| c| c| c| c| c| c|} \hline 
& Quantity & Rank & Publications & Citations & H-Index &  AMS  & Year of PhD \\ \hline
Centroid 1 & 35 &2.571 & 36.686 & 445.400 & 11.029 & .229 & 1999.143 \\ \hline
Centroid 2 & 11 & 3 & 76.455 & 1471.818 & 19.909 & .455 & 1986.636 \\ \hline
Centroid 3 & 9 & 3 & 119.667& 3636.889 & 30.778 & .556 & 1980.667 \\ \hline
Centroid 4 & 3 & 3 & 187.667 & 9024&  38.667& 1 & 1977.333 \\ \hline
\end{tabular}

\textbf{Dartmouth}, \textit{Parameters:} Cosine.

\begin{tabular}{|c ||c| c| c| c| c| c| c|} \hline 
& Quantity & Rank & Publications & Citations & H-Index &  AMS  & Year of PhD \\ \hline
Centroid 1 & 14 & 2.143 & 20.857 & 106.286 & 5.500 & 0 & 1997.929 \\ \hline
Centroid 2 & 6 & 3 & 46.167 & 528.167 & 11.333 & 0 & 1991.667 \\ \hline
Centroid 3 & 3 & 3 & 92.333 & 1211 & 18 & .333 & 1977.666 \\ \hline
\end{tabular}

\textbf{Florida}, \textit{Parameters:} Cosine.

\begin{tabular}{|c ||c| c| c| c| c| c| c|} \hline 
& Quantity & Rank & Publications & Citations & H-Index &  AMS  & Year of PhD \\ \hline
Centroid 1 & 18  & 2 & 21.056 & 64.111 & 4.444 & 0 & 2001.778  \\ \hline
Centroid 2 &  10 &  2.400 & 43.200 & 292.200 & 9.400 & 0 & 1991.700  \\ \hline
Centroid 3 & 9 & 2.889 & 80.222 & 527.333 & 12.222 & .111 & 1983.556 \\ \hline
Centroid 4 & 7 & 3.429 & 98.857 &1357.571 & 16.857 & .143 & 1978.714 \\ \hline
\end{tabular}

\textbf{Harvard}, \textit{Parameters:} Cosine.

\begin{tabular}{|c ||c| c| c| c| c| c| c|} \hline 
& Quantity & Rank & Publications & Citations & H-Index &  AMS  & Year of PhD \\ \hline
Centroid 1 & 10 &3 & 56.700& 1066.900 & 16.900& .300& 1990.600 \\ \hline
Centroid 2 & 7 & 3 & 106.429 & 3070.571 & 30.286 & .571 & 1977.286 \\ \hline
Centroid 3 & 2 & 3 & 313 & 11618 & 46.500 & .500& 1974.500 \\ \hline
\end{tabular}

\textbf{Michigan}, \textit{Parameters:} Cosine.

\begin{tabular}{|c ||c| c| c| c| c| c| c|} \hline 
& Quantity & Rank & Publications & Citations & H-Index &  AMS  & Year of PhD \\ \hline
Centroid 1 & 27 & 2.516 & 23.704 & 162.444 & 6.593 & 0.111 & 1998.111  \\ \hline
Centroid 2 & 12 & 3 & 87.917 & 1313.167 & 18.500 & .500 & 1989.917 \\ \hline
Centroid 3 &  17 & 3.059 & 59.588 & 744.824 & 13.706 & .412 & 1990.412  \\ \hline
Centroid 4 & 6 & 3.667 & 109.333 & 4212 & 27.667 & .833 & 1970 \\ \hline
\end{tabular}

\textbf{MIT}, \textit{Parameters:} Cosine.

\begin{tabular}{|c ||c| c| c| c| c| c| c|} \hline 
& Quantity & Rank & Publications & Citations & H-Index &  AMS  & Year of PhD \\ \hline
Centroid 1 & 18 & 2.056 & 20.167 & 113.222 & 6.167 & .222 & 2004.722  \\ \hline
Centroid 2 & 12 & 2.667 & 36.250 & 503.500 & 12.583 & .250 & 1999.250 \\ \hline
Centroid 3 & 7 & 3 & 182.571 & 5692.857 & 35.714 & .857 & 1972.571  \\ \hline
Centroid 4 & 16 &  3.125 & 80.563 & 1840.938 & 22.250 & .563 & 1993.063 \\ \hline
\end{tabular}

\textbf{Penn}, \textit{Parameters:} Mixed.

\begin{tabular}{|c ||c| c| c| c| c| c| c|} \hline 
& Quantity & Rank & Publications & Citations & H-Index &  AMS  & Year of PhD \\ \hline
Centroid 1 & 6 & 2.000 & 21.000 & 87.333 & 5.333 & .167 & 2006.167  \\ \hline
Centroid 2 & 3 & 3 & 33.667 & 309.333 & 8.667 & .333 & 1993 \\ \hline
Centroid 3 & 6 & 3 & 75.333 & 1289.167 & 18.167 & .500 & 1981.833  \\ \hline
Centroid 4 &  10 & 3.100 & 67 & 664.300 & 14.100 & .500 & 1982.900\\ \hline
\end{tabular}

\textbf{Princeton}, \textit{Parameters:} Cosine.

\begin{tabular}{|c ||c| c| c| c| c| c| c|} \hline 
& Quantity & Rank & Publications & Citations & H-Index &  AMS  & Year of PhD \\ \hline
Centroid 1 &  20 &  1.850 & 18.450 & 274.200 & 7.60 & .250 & 2009.550 \\ \hline
Centroid 2 & 9  &3 & 69.889 & 1732.222 & 21.889 & .667 & 1986.333\\ \hline
Centroid 3 & 13 & 3 & 160.769 & 5240.231 & 35.692 & .615 & 1978.846 \\ \hline
\end{tabular}

\textbf{Rutgers}, \textit{Parameters:} Cosine.

\begin{tabular}{|c ||c| c| c| c| c| c| c|} \hline 
& Quantity & Rank & Publications & Citations & H-Index &  AMS  & Year of PhD \\ \hline
Centroid 1 & 22 & 2.318 &  27.864 &  158.364 & 6.818 &  .273 & 2000.091 \\ \hline
Centroid 2 & 21 & 3.476 & 68.190 & 757.905 & 14.238 & .571 & 1985.762  \\ \hline
Centroid 3 & 7 & 3.714 & 106.571 & 1713.714 & 21.143 &  .857 & 1985.143  \\ \hline
Centroid 4 & 9 & 4 & 159.667 & 3247.556 & 27.222 & 1 & 1972.778  \\ \hline
\end{tabular}

\textbf{UCLA}, \textit{Parameters:} Cosine.

\begin{tabular}{|c ||c| c| c| c| c| c| c|} \hline 
	& Quantity & Rank & Publications & Citations & H-Index &  AMS  & Year of PhD \\ \hline
	Centroid 1 & 23  & 2.478 & 24.783 & 166.130 & 6.565 &  0.127 & 2002.087 \\ \hline
	Centroid 2 & 33 & 2.970 &  70.455 & 1333.212& 17.393 &  0.455& 1989.848 \\ \hline
	Centroid 3 & 2 &  3 & 299.5 & 15861.5 & 58.5 & 1 & 1981 \\ \hline
\end{tabular}

\section{Summary}
In this paper, the exploratory analysis of the math faculty data is conducted and multiple mathematical and statistical methods are used to predict the ranks and AMS fellow status of a math faculty member from other independent variables such as the number of publications and the number of citations. There is a strong demonstration of statistical correlation of the properties examined within the groups, and even with the simpler methods employed, there seems to be much promising potential for the development of an automatic promotion algorithm. For public universities in the United States, salary is listed online and is an additional parameter that may be valuable to predict. We encourage future researchers to make use of the data we have collected and/or additional data and experiment with more refined methods, and academic departments to consider developing and implementing algorithmic promotion methods. 

\section*{Acknowledgement}
We are thankful to Tong Cheng, Terence Coelho, Quentin Dubroff, Joe Olsen and Jason Saied for their contributions in the Experimental Mathematics (Spring 2019) class project \cite{DAMF} at Rutgers, which was the inspiration for this paper.

\bigskip
\hrule
\bigskip
Contact information of the authors:

\{vc362, di110, yao, zeilberger, tz188\} at math dot rutgers dot edu

\end{document}